\documentclass[11pt]{article}
\usepackage{amsmath,graphicx,amsthm,amssymb,algorithm,algorithmic}
\usepackage[
            body={15.2true cm,24true cm},
            headheight=-1.0true cm
            ]{geometry}
\newtheorem{theorem}{Theorem}

\title{On Convergence of the Inexact Rayleigh Quotient Iteration with the
Lanczos Method Used for Solving Linear Systems\footnote{Supported
by National Basic Research Program of China 2011CB302400 and the
National Science Foundation of China (No. 11071140).}}

\author{Zhongxiao Jia\thanks{Department of Mathematical
Sciences, Tsinghua University, Beijing 100084, People's Republic
of China, jiazx@tsinghua.edu.cn}}
\date{}

\begin{document}
\maketitle


\begin{abstract}
For the Hermitian inexact Rayleigh quotient iteration (RQI), the
author has established new local general convergence results, independent
of iterative solvers for inner linear systems. The theory shows that
the method locally converges quadratically under a new condition,
called the uniform positiveness condition. In this paper we
first consider the local convergence of the inexact RQI with the unpreconditioned
Lanczos method for the linear systems. Some attractive properties are
derived for the residuals, whose norms are $\xi_{k+1}$'s,
of the linear systems obtained by the Lanczos
method. Based on them and the new general
convergence results, we make a refined analysis and establish new local
convergence results. It is proved that the inexact RQI with
Lanczos converges quadratically provided that $\xi_{k+1}\leq\xi$ with a
constant $\xi\geq 1$. The method is guaranteed to converge
linearly provided that $\xi_{k+1}$ is bounded by a small multiple of the
reciprocal of the residual norm $\|r_k\|$ of the current approximate eigenpair.
The results are fundamentally different from the existing
convergence results that always require $\xi_{k+1}<1$,
and they have a strong impact on effective
implementations of the method. We extend the new theory to
the inexact RQI with a tuned preconditioned Lanczos for the
linear systems. Based on the new theory,
we can design practical criteria to control $\xi_{k+1}$ to
achieve quadratic convergence and implement
the method more effectively than ever before.
Numerical experiments confirm our theory.
\bigskip


\bigskip


\end{abstract}

\section{Introduction} \label{SecIntro}

We consider the problem of computing an eigenvalue $\lambda$ and the
associated eigenvector $x$ of a large and possibly sparse Hermitian
matrix $A\in \mathbb{C}^{n\times n}$, i.e.,
\begin{equation}
Ax_i=\lambda_i x_i, \ i=1,2,\ldots,n,
\end{equation}
where $(\lambda_i,x_i), \ i=1,2,\ldots,n$ are the eigenpairs
with $\|x_i\|=1$ in the 2-norm.
Throughout the paper, we are interested in the eigenvalue
$\lambda_1$ closest to a target $\sigma$
and its corresponding eigenvector $x_1$ in the sense that
\begin{equation}
|\lambda_1-\sigma| <|\lambda_2-\sigma| \le \cdots \le
|\lambda_n-\sigma|. \label{sigma}
\end{equation}
Suppose that $\sigma$ is between $\lambda_1$ and $\lambda_2$.
Then (\ref{sigma}) means
\begin{equation}
|\lambda_1-\sigma|<\frac{1}{2}|\lambda_1-\lambda_2|.\label{gap}
\end{equation}
There are a number of methods for solving this kind of problem, such as inverse
iteration \cite{ParlettSEP}, Rayleigh quotient iteration (RQI)
\cite{ParlettSEP}, the Lanczos method and its shift-invert variant
\cite{ParlettSEP}, the Davidson method and the Jacobi--Davidson
method \cite{stewart,vorst}. However, except the standard Lanczos
method, these methods and shift-invert Lanczos involve the
solution of a possibly ill-conditioned linear system at each
iteration. This is generally very difficult and even impractical by
a direct solver since a factorization of a shifted $A$ may be too
expensive. So one generally resorts to iterative solvers for
the linear systems, called inner iterations. We call updates of
approximate eigenpairs outer iterations. A combination of them
yields an inner-outer iterative eigensolver, also called
an inexact eigensolver.

For $A$ Hermitian, general local convergence theory on the inexact RQI
can be found in Smit~\cite{Smit}, van den Eshof~\cite{EshofJD}, Simoncini
and Eld$\acute{e}$n~\cite{SimonciniRQI}, Berns-M\"uller and Spence~\cite{mgs06}
and Freitag and Spence~\cite{freitagspence,freitag08b}.
Berns-M\"uller and Spence~\cite{MullerVariableShift}
have extended the theory of \cite{mgs06} to the case that $A$ is non-Hermitian.
For $A$ Hermitian and non-Hermitian,
Xue and Elman~\cite{xueelman} and Xue and Syzld~\cite{xueszyld} have
further analyzed the local convergence of the inexact RQI and provided new
insights into tuning a preconditioner, which is for efficient Krylov inner solves.
The idea of tuning a preconditioner was initially proposed
in \cite{SimonciniRQI} and then extended in \cite{mgs06,MullerVariableShift}
and improved in \cite{freitagspence,freitag08b}.
Let $\|r_k\|$ be the residual norm of the approximate eigenpair at outer
iteration $k$. The mentioned papers except \cite{xueszyld} have proved
that the inexact RQI converges cubically
if $\xi_{k+1}=O(\|r_k\|)$ and quadratically if $\xi_{k+1}\leq\xi<1$ with a
constant $\xi$ not near one. In  \cite{xueszyld}, Xue and Syzld have given
a new local convergence analysis, showing that the inexact RQI may
demonstrate cubic and quadratic asymptotic convergence rates, respectively,
for Hermitian and non-Hermitian problems, if the shifted linear
systems are solved by a Krylov subspace method with a tuned preconditioner
to a reasonably small fixed tolerance $\xi_{k+1}\leq\xi<1$. However, their
claims hold only under a crucial assumption that the factor $\alpha_m$ in their
main result (3.2) must be bounded by a moderate constant. However, it
is not clear when this assumption is satisfied. They have given some qualitative
but rough and non-rigorous arguments on $\alpha_m$, trying to show that $\alpha_m$
is bounded.
In fact, the size of $\alpha_m$ is closely related to that of $\xi_{k+1}$, as was
implicit from their arguments. Furthermore, their arguments
implicitly indicate that $\alpha_m$ is guaranteed to
be moderate only when $\xi_{k+1}$ is sufficiently small, and
there is no evidence that a reasonably
small fixed $\xi_{k+1}$ is enough. So (3.2) in Theorem 3.1 of
\cite{xueszyld} may not mean the cubic asymptotic convergence if the shifted linear
systems are solved by a tuned preconditioned Krylov subspace method
to a reasonably small fixed tolerance $\xi_{k+1}\leq\xi<1$. Particularly,
if the Lanczos method is used for solving shifted inner linear systems,
it is typical that $\alpha_m$ can be very large and $\xi_{k+1}$ may be bigger
than one for some inner iterations when the linear
systems are indefinite.\footnote{I have communicated with Dr Fei Xue,
one author of \cite{xueszyld}. He has agreed that (i) their
arguments on $\alpha_m$ are non-rigorous and more experimental and (ii)
the condition of the cubic convergence result is not yet clear
and how small $\xi_{k+1}$ should be is being under consideration.}

Simoncini and Eld$\acute{e}$n~\cite{SimonciniRQI} are the first to observe
that poor approximate solutions of the linear systems may be much improved
approximations to the desired eigenvector $x_1$. In our notation, their
observation qualitatively means that $\xi_{k+1}$ can be near
one when MINRES is used. More generally, a remarkable feature for the
convergence of the inexact RQI with MINRES (or its mathematically
equivalent version the conjugate residual method) is that $\xi_{k+1}$
is allowed not small. Xue and Elman~\cite{xueelman} have given a {\em qualitative}
mathematical justification on this phenomenon. We must point out that their results
have also shown clearly that although the convergence of the inexact RQI with
MINES allows $\xi_{k+1}$ near one, $\xi_{k+1}$ too near one is prohibited.
Quantitatively, however, how near one $\xi_{k+1}$ is allowed is unknown.
The observation by Simoncini and Eld$\acute{e}$n and the work followed
may have led to some serious misunderstanding or wrong impressions that the size of
$\xi_{k+1}$ plays {\em no} role in the convergence of the inexact
RQI. In fact, this is completely wrong, and the size of
$\xi_{k+1}$ {\em does} play a {\em crucial} role in the convergence and convergence
rates of the inexact RQI. It can be easily observed that in literature
all the convergence conditions involve
$\xi_{k+1}$ either explicitly or implicitly, so do stopping criteria.
They all assume $\xi_{k+1}<1$ first and then impose further restrictions to $\xi_{k+1}$
either explicitly or implicitly for achieving a desired convergence rate.
This is easily understandable since an approximate solution corresponds to a {\em unique}
$\xi_{k+1}$. Mathematically, any result on approximate solutions
can be definitely reflected by their residuals, and $\xi_{k+1}$
plays a role whenever an approximate solution
of the linear system is invoked. Ignoring or overlooking the role of
$\xi_{k+1}$ is obviously incorrect.
For example, Xue and Elman~\cite{xueelman} have adapted the stopping criteria
from \cite{SimonciniRQI}.
Seemingly, the criterion does not involve $\xi_{k+1}$ explicitly but it requires
that $p_m(\lambda_1-\theta_k)$ be reasonably small,
say $10^{-3}$ to $10^{-2}$, where $p_m(z)$ is the associated residual polynomial
of $m$-step MINRES for the linear system and
$\theta_k$ is the Rayleigh quotient of the inexact RQI at step $k$. However,
the size of $p_m(\lambda_1-\theta_k)$ is very closely related to $\xi_{k+1}$.
Actually, it is easy to justify from the theory of MINRES that a small
$|p_m(\lambda_1-\theta_k)|$ must mean a small $\xi_{k+1}$, and vice versa.
Therefore, although $\xi_{k+1}$ is not necessarily small, its {\em size} does
play a key role {\em in its own intrinsic way} in determining the cubic,
quadratic and linear convergence rates of the inexact RQI.
It is nontrivial and delicate for us to find out the correct
quantitative role of $\xi_{k+1}$ and provide more insights into the convergence
rates of the inexact RQI.

Before \cite{xueszyld} was available in a technical report form in September 2010,
the author of this paper in June 2009
gave a refined analysis on the inexact RQI with MINRES used for
solving the linear systems, and the latest third version was
available in June 2010 \cite{jia09}. It is proved that the RQI with MINRES
generally achieves the cubic asymptotic convergence whenever
$\xi_{k+1}$ is not near one. The quadratic asymptotic convergence
quantitatively requires $\xi_{k+1}=1-O(\|r_k\|)$ (here we express it qualitatively),
which is increasingly nearer to one as outer iterations proceeds and
is much more relaxed than the existing ones
in literature. Also, the linear convergence conditions are given
that quantitatively require only $\xi_{k+1}=1-O(\|r_k\|^2)$, nearer to
one than for the quadratic convergence.
For the inexact RQI, to the author's best
knowledge, there has been no result available on linear convergence.
Besides \cite{jia09}, numerical experiments
in \cite{mgs06,xueszyld} have also confirmed such cubic asymptotic
convergence.

As we have noticed, the basic condition $\xi_{k+1}\leq\xi<1$ is first assumed
in all the above mentioned papers except \cite{jia09}. This condition,
though seemingly natural and necessary, may miss something essential
and prevents us from getting better and more insightful results on
the convergence of the inexact RQI.
The author in \cite{jia09} has revisited the convergence of the inexact
RQI independent of iterative solvers and presents new general local
convergence results. It is proved that the inexact RQI
converges quadratically under a so-called uniform positiveness
condition, which retains more information on inner solves
and is fundamentally different from and weaker
than the condition $\xi_{k+1}\leq\xi<1$.
Although several results have been established for the
inexact RQI with Lanczos in literature, one treats
the residuals obtained by the Lanczos method as general
ones and simply takes their norms in convergence analysis. Therefore,
residual directions of inner iterations have not been reasonably
exploited, and {\em fundamental} effects of residual directions
on convergence have been completely overlooked.
Based on the new general convergence results in \cite{jia09},
we first establish some attractive properties of the residuals
obtained by the Lanczos method for the linear systems. By fully
exploiting them, we then make a novel analysis and
derive a number of insightful results that are not only
stronger than but also fundamentally different
from the ones available in literature.

We should stress an important fact that each shifted inner linear system
involved in the inexact RQI is typically Hermitian
indefinite. For indefinite problems, it is well known that
the Lanczos method typically behaves irregularly, that is,
residual norms $\xi_{k+1}$'s of the approximate solutions
of the linear system can be large and even infinite
for some Lanczos steps. For the inexact RQI with the
unpreconditioned Lanczos, the most remarkable results
we will prove are that the inexact RQI
with Lanczos asymptotically converges quadratically
provided that $\xi_{k+1}\leq\xi$ with a constant $\xi$ that is allowed to
be {\em bigger} than one and the method
is guaranteed to converge linearly provided that $\xi_{k+1}$ is
bounded by not exceeding a moderate multiple of $\frac{1}{\|r_k\|}$,
which means that we allow $\xi_{k+1}\gg 1$.
The results have a strong impact on effective implementations of the
method. Based the new theory, we design new stopping criteria for
inner solves. To achieve the quadratic convergence, compared with prevailing
implementations of the method, our new implementation can
save the computational cost of solving the
linear systems very significantly.
Numerical experiments demonstrate that the new implementation is
twice to four times and even more as fast as
the prevailing implementations and the method converges
smoothly and quickly for some problems
even when $\xi_{k+1}$ is up to $10^4\sim 10^7$.

As byproducts, similar to that done in \cite{jia09,SimonciniRQI},
we establish lower bounds on the norms of approximate solutions
$w_{k+1}$ of the linear systems obtained by the unpeconditioned Lanczos.
We show that $\|w_{k+1}\|$ is always
$O(\frac{1}{\|r_k\|^2})$ no matter the inexact RQI with Lanczos
converges cubically or quadratically. Therefore, it is distinctive that
$\|w_{k+1}\|$ itself obtained by Lanczos cannot reveal the convergence
behavior of the inexact RQI and cannot be used to design stopping
criteria for inner iterations. Making use of these bounds, we
present a simpler but weaker quadratic convergence result.
As a global result, similar to that for
the inexact RQI with MINRES \cite{SimonciniRQI} where it is
shown that outer residual norms $\|r_k\|$ do not decrease
monotonically any more for an {\em arbitrary} starting vector,
we derive a relationship between $\|r_k\|$ and $\|r_{k+1}\|$,
starting with an arbitrary vector instead of a reasonably good one.
We will see that, unlike the exact RQI, the inexact RQI with Lanczos
loses the residual monotonic decreasing property for an arbitrary
starting vector. Therefore, {\em for the convergence of the inexact RQI,
it is only meaningful to speak of local rather than global convergence.
That is, under the assumption that the current approximation
has already a reasonable accuracy, one investigates how the
next approximation better approaches the desired eigenvector.
By convergence (rate), we always mean asymptotic convergence (rate)}.

We also extend our theory to the inexact RQI with a tuned
preconditioned Lanczos method. This is a nontrivial task. We
show that our main results in the unpreconditioned case can be extended to the
tuned preconditioned case.

The paper is organized as follows. In Section \ref{SecIRQI}, we
review the inexact RQI and the new general convergence theory of
\cite{jia09} on the inexact RQI. In
Section~\ref{seclanczos}, we present convergence results on the
inexact RQI with the unpreconditioned Lanczos Lanczos for solving
inner linear systems. In
Section~\ref{precondit}, we extend the theory to the inexact RQI
with a tuned preconditioned Lanczos method for solving inner linear
systems. We perform numerical experiments to
confirm our results in Section~\ref{testlanczos}.
Finally, we end up with some concluding remarks in Section~\ref{conc}.

Throughout the paper, denote by the superscript
* the conjugate transpose of a matrix or vector,
by $\|\cdot\|$ the vector 2-norm and
the matrix spectral norm, and by $\lambda_{\min},\lambda_{\max}$
the algebraically smallest and largest eigenvalues of $A$, respectively.

\section{The inexact RQI and general convergence theory} \label{SecIRQI}

RQI is a famous iterative algorithm and its locally cubic
convergence for Hermitian problems is very attractive
\cite{ParlettSEP}. It plays a crucial role in some practical
effective algorithms, e.g., the QR algorithm,
\cite{GolubMC,ParlettSEP}. Assume that the unit length $u_k$ is
{\em already} a reasonably good approximation
to $x_1$. Then the Rayleigh quotient $\theta_k = u^*_k A u_k$ is a
good approximation to $\lambda_1$ too. RQI \cite{GolubMC,ParlettSEP}
computes a new approximation $u_{k+1}$ to $x_1$ by solving the shifted
inner linear system
\begin{equation} \label{EqERQILinearEquation1}
( A - \theta_k I ) w = u_k
\end{equation}
for $w_{k+1}$ and updating $u_{k+1}=w_{k+1}/\|w_{k+1}\|$ and
iterates until convergence. It is known
\cite{mgs06,NotayRQI,ParlettSEP} that if
$$
|\lambda_1-\theta_0|
<\frac{1}{2}\min_{j=2,3,\ldots,n}|\lambda_1-\lambda_j|
$$
then RQI asymptotically converges to $\lambda_1$ and $x_1$
cubically. So we can assume that the eigenvalues of $A$ are ordered
as
\begin{equation}
|\lambda_1-\theta_k| <|\lambda_2-\theta_k| \le \cdots \le
|\lambda_n-\theta_k| \mbox{ for all $k$}. \label{order}
\end{equation}
With this ordering and noting that
$\lambda_{\min}\leq\theta_k\leq\lambda_{\max}$, we have
\begin{equation}
|\lambda_1-\theta_k|<\frac{1}{2}|\lambda-\lambda_2|.\label{sep1}
\end{equation}

In the inexact RQI, \eqref{EqERQILinearEquation1} is solved by an iterative solver
and an approximate solution $w_{k+1}$ satisfies
\begin{equation} \label{EqIRQILinearEquation1}
( A - \theta_k I )w_{k+1} = u_k + \xi_{k+1} d_{k+1}, \quad u_{k+1} =w_{k+1}
/ \|w_{k+1} \|
\end{equation}
with $ 0 < \xi_{k+1} \le \xi$, where $\xi_{k+1}d_{k+1} $ with $\|d_{k+1}\|=1$ is the
residual of $(A-\theta_k I)w=u_k$, $d_{k+1}$ is the residual direction
vector and $\xi_{k+1}$ is the {\em relative}
residual norm (inner tolerance) as $\|u_k\|=1$ and may change at
every outer iteration $k$. This process is summarized as Algorithm
1. If $\xi_{k+1}=0$ for all $k$, Algorithm 1 becomes the exact RQI.

\begin{algorithm}
  \caption{The inexact RQI} \label{AlgIRQI}
  \begin{algorithmic}[1]
  \STATE Choose a unit length $ u_0 $, a reasonable approximation to $ x_1$.
  \FOR{$k$ = 0,1, \ldots}
  \STATE $\theta_k =u^*_k A u_k$.
  \STATE Solve $(A-\theta_k I)w=u_k$ for $w_{k+1}$ by an
  iterative solver with
  $$
  \|(A-\theta_k I )
  w_{k+1}-u_k\|=\xi_{k+1}.
  $$
  \STATE $u_{k+1}=w_{k+1}/\|w_{k+1} \|$.
  \STATE If convergence occurs, stop.
  \ENDFOR
  \end{algorithmic}
\end{algorithm}

There are a number of general local quadratic convergence results in, e.g.,
\cite{mgs06,SimonciniRQI,Smit, EshofJD}, which are all obtained by first requiring
$\xi_{k+1}\leq\xi<1$. In \cite{jia09}, new local general
convergence results have been proved under a new condition that is
fundamentally different and can relax $\xi_{k+1}$ very much.
To present the results, we decompose $u_k$ and $d_{k+1}$
into the orthogonal direct sums
\begin{eqnarray}
&u_k = x \, \cos \phi_k + e_k \, \sin \phi_k, \quad e_k \perp x,
\label{EqIRQIDecompositoinOfu_k} \\
&d_{k+1} = x \, \cos \psi_k + f_k \, \sin \psi_k, \quad f_k \perp x
\label{EqIRQIDecompositoinOfd_{k+1}}
\end{eqnarray}
with $\|e_k\|=\|f_k\|=1$ and $\phi_k=\angle(u_k,x)$, $\psi_k=\angle(d_{k+1},x)$.
Here without loss of generality and for brevity of discussions,
we suppose that $\phi_k$ is the acute angle
between $u_k$ and $x_1$. Furthermore, we stress again that
speaking of local convergence analysis naturally means that $\phi_k$ is
already reasonably small, i.e., $\cos\phi_k\approx 1$ and $\sin\phi_k\approx 0$,
and one then investigates how the next $u_{k+1}$ and $\theta_{k+1}$ better
approximate $x_1$ and $\lambda_1$.

Given this, we should remind that
$\cos\psi_k$ is either positive or negative depending on
$d_{k+1}$. Note that \eqref{EqIRQILinearEquation1} can be written as
\begin{eqnarray} \label{EqIRQILinearEquation2}
( A - \theta_k I ) w_{k+1} = ( \cos \phi_k + \xi_{k+1} \, \cos \psi_k )
\, x + ( e_k \, \sin \phi_k + \xi_{k+1} \, f_k \, \sin \psi_k ).
\end{eqnarray}
Inverting $ A - \theta_k I$ gives
\begin{equation} \label{EqIRQIwk}
w_{k+1} = (\lambda_1 - \theta_k )^{-1} ( \cos \phi_k + \xi_{k+1} \, \cos
\psi_k ) \, x + ( A - \theta_k I )^{-1} ( e_k \, \sin \phi_k + \xi_{k+1}
\, f_k \, \sin \psi_k ).
\end{equation}
Define $\|r_k\|=\|(A-\theta_k I)u_k\|$. Then by (\ref{sep1}) we get
$|\lambda_2-\theta_k|>\frac{|\lambda_2-\lambda_1|}{2}$. It is known
from \cite[Theorem 11.7.1]{ParlettSEP} that
\begin{equation}
\frac{\|r_k\|}{\lambda_{\max}-\lambda_{\min}}\leq\sin\phi_k\leq\frac{2\|r_k\|}
{|\lambda_2-\lambda_1|}.\label{parlett}
\end{equation}
We comment that $\lambda_{\max}-\lambda_{\min}$ is the spectrum
spread of $A$ and $|\lambda_2-\lambda_1|$ is the gap or separation of
$\lambda_1$ and the other eigenvalues of $A$.

Throughout the paper, we define
\begin{equation}
\beta=\frac{\lambda_{\max}-\lambda_{\min}}{|\lambda_2-\lambda_1|}.
\label{beta}
\end{equation}

\begin{theorem}{\rm \cite{jia09}} \label{ThmIRQIQuadraticConvergence}
If the uniform positiveness condition
\begin{equation} \label{EqIRQIC}
|\cos \phi_k + \xi_{k+1} \cos \psi_k| \ge c
\end{equation}
is satisfied with a moderate constant $c>0$ independent of $k$,
then
\begin{eqnarray}
\tan \phi_{k+1} &\le &2\beta\frac{\sin \phi_k + \xi_{k+1} \sin \psi_k}
{|\cos \phi_k + \xi_{k+1} \cos \psi_k|}\sin^2\phi_k   \label{bound1}\\
&\le &\frac{2\beta\xi_{k+1}}{c} \sin^2 \phi_k + O( \sin^3 \phi_k ),
\label{EqIRQIQuadraticConvergence}
\end{eqnarray}
that is, the inexact RQI asymptotically converges quadratically provided that
{\rm (\ref{EqIRQIC})} is satisfied and $\xi_{k+1}$ is uniformly bounded by
some moderate constant.
\end{theorem}

It can be found in \cite{jia09} that the proof of (\ref{bound1}) is elementary
and easy to follow. Combining (\ref{bound1}) and (\ref{EqIRQIC}), it
is direct to get (\ref{EqIRQIQuadraticConvergence}).

\begin{theorem}{\rm \cite{jia09}}\label{resbound}
If the uniform positiveness condition {\rm (\ref{EqIRQIC})} holds,
then
\begin{equation}
\|r_{k+1}\|\leq\frac{8\beta^2\xi_{k+1}} {c|\lambda_2-\lambda_1|}\|r_k\|^2
+O(\|r_k\|^3).\label{resgeneral}
\end{equation}
\end{theorem}

We make some comments on the above two theorems.

{\bf Remark 1}. They illustrate that it is the size of $|\cos\phi_k
+ \xi_{k+1}\cos \psi_k|$ other than $\xi_{k+1}\leq\xi<1$ that is critical in
convergence.

{\bf Remark 2}. If $ \xi_{k+1} = 0 $ for all $ k $, then the inexact RQI
reduces to the exact RQI and
Theorems~\ref{ThmIRQIQuadraticConvergence}--\ref{resbound} show the
cubic convergence: $\tan
\phi_{k+1}\leq\frac{2\beta}{\cos\phi_k}\sin^3\phi_k$ and
$\|r_{k+1}\|=O(\|r_k\|^3)$. If the linear systems are solved with
decreasing tolerance $\xi_{k+1}=O(\sin\phi_k)=O(\|r_k\|)$, then we have
the cubic convergence: $\tan \phi_{k+1}=O(\sin^3\phi_k)$ and
$\|r_{k+1}\|=O(\|r_k\|^3)$.

{\bf Remark 3}. If $\cos\psi_k$ is positive, the uniform
positiveness condition holds for any uniformly bounded
$\xi_{k+1}\leq\xi$ with $\xi$ a moderate constant. So we may have
$\xi\geq 1$. If $\cos\psi_k$ is negative, the uniform
positiveness condition $|\cos\phi_k+\xi_{k+1}\cos\psi_k|\geq c$ means
that
$$
\xi_{k+1}\leq\frac{c-\cos\phi_k}{\cos\psi_k}
$$
if $\cos\phi_k+\xi_{k+1}\cos\psi_k\geq c$ and
$$
\xi_{k+1}\geq\frac{c+\cos\phi_k}{-\cos\psi_k}
$$
if $-\cos\phi_k-\xi_{k+1}\cos\psi_k\geq c$.
So the size of $\xi_{k+1}$ critically depends on
that of $\cos\psi_k$, and for a given $c$ we may have $\xi_{k+1}\approx 1$
and even $\xi_{k+1}>1$. Obviously, without the information on $\cos\psi_k$, it
would be impossible to access or estimate $\xi_{k+1}$. As a general
convergence result, however, its significance and importance consist
in that it fully exploits a crucial quantity $\cos\psi_k$ to relax
$\xi_{k+1}$
as much as possible and meanwhile preserves the same convergence rate
of outer iteration. As a result, the condition $\xi_{k+1}\leq\xi<1$ with
constant $\xi$ not near one may be stringent and unnecessary for the
quadratic convergence of the inexact RQI, independent of iterative
solvers for inner linear systems.

{\bf Remark 4}. Keep in mind a basic fact that the bigger $\xi_{k+1}$ is,
the less costly a chosen inner iterative solver is.
The new conditions on $\xi_{k+1}$ derived from the uniform positiveness
condition have a strong impact on effective implementations of the inexact RQI
since we must stop a certain iterative solver, e.g.,
the very popular MINRES method and the Lanczos method for solving
$(A-\theta_k I)w=u_k$ at {\em right} moment. It appears
that $\cos\psi_k$ is critically iterative solver dependent. For Lanczos,
$\cos\psi_k$ has some very attractive properties. Making use of them,
we can precisely determine bounds for $\xi_{k+1}$ in
Section~\ref{seclanczos}, which are much more relaxed than those
in literature. For the convergence of the inexact RQI with MINRES,
we refer to \cite{jia09} for the properties of $\cos\psi_k$ and their effects on
$\xi_{k+1}$.

\section{Convergence of the inexact RQI with the unpreconditioned Lanczos}
\label{seclanczos}

The previous results and discussions are for general purpose,
independent of iterative solvers for $(A-\theta_k I)w=u_k$. Since
we have $\lambda_{\rm min}\leq\theta_k\leq\lambda_{\rm max}$,
the matrix $A-\theta_k I$ is Hermitian indefinite when $\theta_k\not=\lambda_1$.
The Lanczos method is a popular Krylov subspace iterative solver for
Hermitian linear systems \cite{saad}. The method nicely fits
into the inexact RQI.

We briefly review the Lanczos method for solving
(\ref{EqERQILinearEquation1}). At outer iteration $k$, taking the
starting vector $v_1$ to be $u_k$, the $m$-step Lanczos process on
$A-\theta_k I$ can be written as
\begin{equation}
( A - \theta_k I ) V_m = V_mT_m+t_{m+1m}v_{m+1}e_m^*,
\label{lanczosp}
\end{equation}
where the columns of $ V_m =(v_1,\ldots,v_m)$ form an orthonormal
basis of the Krylov subspace $\mathcal{K}_m(A - \theta_k I, u_k )
= \mathcal{K}_m ( A, u_k )$ and $T_m=(t_{ij})=V_m^*(A-\theta_k
I)V_m$ is an $m\times m $ Hermitian tridiagonal matrix
\cite{ParlettSEP,saad}.

The Lanczos method \cite{GolubMC,paige,saad}
is a Galerkin projection method and requires the
residual $\xi_{k+1}d_{k+1}$ to be orthogonal to the search subspace.
With the zero vector as an initial guess to the solution of
$(A-\theta_k I)w=u_k$, the Galerkin condition means
$\xi_{k+1}d_{k+1}\perp {\cal K}_m(A,u_k)$. Specially,
$\xi_{k+1}d_{k+1}\perp u_k\in {\cal K}_m(A,u_k)$.
The method extracts the approximate solution
$w_{k+1}=V_m\hat y$ to $(A-\theta_k I)w=u_k$ from $\mathcal{K}_m
(A, u_k)$, where $\hat y$ is the solution of the Hermitian
tridiagonal linear system $T_my=e_1$ with $e_1$ being the first
coordinate vector of dimension $m$. It is worth noting that we
should naturally take $m>1$; otherwise
we would have $\mathcal{K}_1 ( A, u_k )={\rm span}\{u_k\}$ and
$T_1=0$, so that the Lanczos method would break down and $w_{k+1}$
would not exist. The algorithm SYMMLQ is a very effective implementation of the
Lanczos method \cite{GolubMC,paige75}.

For Hermitian positive definite linear systems,
the Lanczos method is mathematically equivalent to the conjugate
gradient method and has the optimality that the error of the
approximate solution is minimal with respect to the energy norm over
the given Krylov subspace \cite{saad}. For Hermitian indefinite
linear systems, the method does not have any kind of optimality.
For our case, the linear system $(A-\theta_kI)w=u_k$ is not only
indefinite but also increasingly ill conditioned as
$\theta_k\rightarrow\lambda_1$. The indefinite
system is typically ill conditioned and can be (nearly) singular,
so that the Lanczos method may converge slowly and
irregularly and $\xi_{k+1}$ can be typically big for $m$ small. For more
details, we refer to \cite{paige,saad}.

We will present convergence results on the inexact
RQI with Lanczos. First of all, we establish the following results,
which will play a key role in the later analysis.

\begin{theorem}\label{lanczos}
It holds that
\begin{equation}
|\cos\psi_k|\leq\tan\phi_k \label{lanczos1}
\end{equation}
and asymptotically
\begin{equation}
\sin\psi_k\geq 1-\frac{1}{2}\sin^2\phi_k\label{lanczos2}
\end{equation}
by ignoring the higher order term $O(\sin^4\phi_k)$.
\end{theorem}

\begin{proof}
Recall that its residual $\xi_{k+1}d_{k+1}$ obtained by the Lanczos method
satisfies $\xi_{k+1}d_{k+1}\perp {\cal K}_m(A,u_k)$. So, we specially have
$d_{k+1}\perp u_k$. Therefore, from \eqref{EqIRQIDecompositoinOfu_k}
and \eqref{EqIRQIDecompositoinOfd_{k+1}} we get
$$
\cos\phi_k\cos\psi_k+e_k^*f_k\sin\phi_k\sin\psi_k=0,
$$
which means
$$
\frac{|\cos\psi_k|}{\sin\psi_k}=|e_k^*f_k\tan\phi_k|\leq\tan\phi_k.
$$
From the above it follows that (\ref{lanczos1}) holds. By the Taylor
expansion, from (\ref{lanczos2}) we get
$$
\sin\psi_k=\sqrt{1-\cos^2\psi_k}=1-\frac{1}{2}|\cos^2\psi_k|+
O(\cos^4\psi_k)\geq
1-\frac{1}{2}\tan^2\phi_k=1-\frac{1}{2}\sin^2\phi_k
$$
by dropping the higher order term $O(\sin^4\phi_k)$.
\end{proof}

Combining this theorem with (\ref{bound1}) of
Theorem~\ref{ThmIRQIQuadraticConvergence}, we can establish one of
our main results for the inexact RQI with Lanczos.

\begin{theorem}\label{convlanczos}
Let $\xi$ be a constant such that $\xi_{k+1}\leq\xi$ and satisfy
\begin{equation}
\xi\sin\phi_k\leq\alpha<1 \label{suffcond1}
\end{equation}
with $\alpha$ a given constant not near one.
Then the uniform positiveness condition
{\rm (\ref{EqIRQIC})} holds and the inexact RQI with Lanczos
asymptotically converges quadratically:
\begin{equation}
\tan\phi_{k+1}\leq \frac{2\beta \xi}{1-\alpha}\sin^2\phi_k.
\label{lanczosquadra2}
\end{equation}
It asymptotically converges cubically if
\begin{equation}
\xi_{k+1}=O(\sin\phi_k); \label{lanczoscubic}
\end{equation}
it converges at linear factor $\gamma<1$:
\begin{equation}
\tan\phi_{k+1}\leq\gamma\sin\phi_k \label{linear}
\end{equation}
if
\begin{equation}
\xi_{k+1}\leq\frac{\gamma-2\beta\sin^2\phi_k}{(2\beta+\gamma)\sin\phi_k}
<\frac{\gamma}{(2\beta+\gamma)\sin\phi_k}. \label{linconv}
\end{equation}
\end{theorem}

\begin{proof}
By (\ref{lanczos2}), we get
$$
\sin\phi_k+\xi_{k+1}\sin\psi_k\leq\sin\phi_k+\xi_{k+1}\leq\sin\phi_k+\xi.
$$
On the other hand, we get from (\ref{lanczos1})
\begin{eqnarray}
|\cos\phi_k+\xi_{k+1}\cos\psi_k| &\geq&|\cos\phi_k-\xi_{k+1}|\cos\psi_k||\nonumber\\
&\geq&
|1-\frac{1}{2}\sin^2\phi_k+O(\sin^4\phi_k)-\xi_{k+1} \tan\phi_k| \nonumber\\
&=&|1-\xi_{k+1}\sin\phi_k+O(\sin^2\phi_k)| \nonumber\\
&\geq&1-\xi\sin\phi_k \nonumber\\
&\geq&1-\alpha, \label{denolanczos}
\end{eqnarray}
by dropping the higher order term $O(\sin^2\phi_k)$.
So the uniform positiveness condition (\ref{EqIRQIC}) holds with $c=1-\alpha$.
We then derive from
(\ref{bound1}) that
\begin{eqnarray}
\tan\phi_{k+1}&\leq&2\beta
\frac{\sin\phi_k+\xi_{k+1}\sin\psi_k}{|\cos\phi_k+\xi_{k+1}\cos\psi_k|}\sin^2\phi_k\nonumber\\
&\leq&2\beta\frac{(\sin\phi_k+\xi_{k+1})}
{|1-\xi_{k+1}\sin\phi_k|}\sin^2\phi_k \label{lanczosquadra}\\
&\leq&2\beta\frac{(\sin\phi_k+\xi)}
{1-\alpha}\sin^2\phi_k \nonumber  \\
&\leq&2\beta\frac{\xi}
{1-\alpha}\sin^2\phi_k+O(\sin^3\phi_k) \nonumber
\end{eqnarray}
which is just (\ref{lanczosquadra2}) by ignoring $O(\sin^3\phi_k)$.

The cubic asymptotic convergence is direct from (\ref{lanczosquadra}) if
$\xi_{k+1}=O(\sin\phi_k)$.

It follows from (\ref{lanczosquadra}) that the inexact RQI with
Lanczos converges at linear factor $\gamma$ at least if for all
$k$ it holds that $\xi_{k+1}\sin\phi_k<1$ and
$$
2\beta\frac{(\sin\phi_k+\xi_{k+1})\sin\phi_k}{1-\xi_{k+1}\sin\phi_k}\leq\gamma<1,
$$
from which we get condition (\ref{linconv}) by manipulation.
\end{proof}

Theorem~\ref{convlanczos} presents the conditions on cubic,
quadratic and linear convergence in terms of an a priori uncomputable
$\sin\phi_k$. We next give their alternatives in terms of the
computable $\|r_k\|$, so that they are of practical value as much as
possible and can be used to control the inner tolerance to achieve a
desired convergence rate.

\begin{theorem}\label{conds}
Let $\xi$ be a constant such that $\xi_{k+1}\leq\xi$ for all $k$ and satisfy
\begin{equation}
\frac{2\xi\|r_k\|}{|\lambda_2-\lambda_1|}\leq\alpha<1 \label{suffcond2}
\end{equation}
with $\alpha$ a given constant not near one.
The uniform positiveness condition holds
and the inexact RQI with Lanczos
asymptotically converges quadratically:
\begin{eqnarray}
\|r_{k+1}\|&\leq&\frac{8\beta^2(2\|r_k\|+\xi_{k+1}|\lambda_2-\lambda_1|)}
{|\lambda_2-\lambda|(|\lambda_2-\lambda_1|-2\xi_{k+1}\|r_k\|)}
\|r_k\|^2, \label{resquadra}\\
&\leq&\frac{8\beta^2(2\|r_k\|+\xi|\lambda_2-\lambda_1|)}
{(\lambda_2-\lambda_1)^2(1-\alpha)} \|r_k\|^2.\label{resquadra2}
\end{eqnarray}
It asymptotically converges cubically if
\begin{equation}
\xi_{k+1}=O(\|r_k\|);\label{cubiclanczos}
\end{equation}
it converges at linear factor $\gamma<1$:
\begin{equation}
\|r_{k+1}\|\leq\gamma\|r_k\| \label{linearr2}
\end{equation}
if
\begin{equation}
\xi_{k+1}\leq\frac{\gamma(\lambda_2-\lambda_1)^2-16\beta^2\|r_k\|^2}
{2|\lambda_2-\lambda_1|(4\beta^2+\gamma)\|r_k\|} < \frac{\gamma
|\lambda_2-\lambda_1|}{(8\beta^2+2\gamma)\|r_k\|}. \label{linconv2}
\end{equation}
\end{theorem}

\begin{proof}
Making use of (\ref{parlett}) gives
$$
\frac{\|r_{k+1}\|}{\lambda_{\max}-\lambda_{\min}}\leq\sin\phi_{k+1}\leq\tan\phi_{k+1},
$$
$$
1-\xi_{k+1}\sin\phi_k\geq
1-\xi_{k+1}\frac{2\|r_k\|}{|\lambda_2-\lambda_1|}\geq
1-\frac{2\xi\|r_k\|}{|\lambda_2-\lambda_1|}\geq 1-\alpha>0
$$
and
$$
\sin\phi_k+\xi_{k+1}\leq\frac{2\|r_k\|}{|\lambda_2-\lambda_1|}+\xi_{k+1}.
$$
Substituting the above relations into (\ref{lanczosquadra}) and
(\ref{lanczosquadra2}) establishes (\ref{resquadra}) and
(\ref{resquadra2}), respectively. It is clear from
(\ref{resquadra}) that the inexact RQI with Lanczos asymptotically converges
cubically once $\xi_{k+1}=O(\|r_k\|)$.

In order to make $\|r_k\|$ linearly converge to zero
monotonically, from (\ref{resquadra}) we simply set
$$
\frac{8(\lambda_n-\lambda)^2}{|\lambda_2-\lambda_1|^3}
\frac{2\|r_k\|+\xi_{k+1}|\lambda_2-\lambda_1|}{|\lambda_2-\lambda_1|-2\xi_{k+1}\|r_k\|}
\|r_k\|\leq\gamma<1.
$$
Solving it for $\xi_{k+1}$ gives
$$
\xi_{k+1}\leq\frac{\gamma(\lambda_2-\lambda_1)^2-16\beta^2\|r_k\|^2}
{2|\lambda_2-\lambda_1|(4\beta^2+\gamma)\|r_k\|} <
\frac{\gamma|\lambda_2-\lambda_1|} {(8\beta^2+2\gamma)\|r_k\|}.
$$
\end{proof}

We make some comments on Theorems~\ref{convlanczos}--\ref{conds}.
\smallskip

{\bf Remark 1}. The quadratic asymptotic convergence condition (\ref{suffcond1})
indicates that $\xi>1$ is allowed as
$\xi\leq\frac{\alpha}{\sin\phi_k}$ and $\sin\phi_k$ is supposed to small.
For a given reasonably good starting vector $u_0$, if both the global
convergence and quadratic asymptotic convergence are required,
then only (\ref{suffcond1}) may not be sufficient. It is seen from
(\ref{lanczosquadra2}) that if $\xi$ satisfies
 \begin{equation}
\tan\phi_{1}\leq \frac{2\beta \xi}{1-\alpha}\sin^2\phi_0<\tan\phi_0
\label{xicond}
\end{equation}
then $\tan\phi_k$ decreases from the beginning of outer iteration.
From (\ref{xicond}) we find
$$
\xi<\frac{1-\alpha}{\beta\sin 2\phi_0}.
$$
Combining it with the requirement $\xi\leq\frac{\alpha}{\sin\phi_0}$, we get
\begin{equation}
\xi<\min\{\frac{\alpha}{\sin\phi_0},\frac{1-\alpha}{\beta\sin 2\phi_0}\}.
\label{xichoice}
\end{equation}
With such $\xi$, the inexact RQI with Lanczos achieves the quadratic asymptotic
convergence.
The upper bound depends on $\beta$,
which measures the conditioning of $x_1$. We see that $\xi$ may be bigger than one
as $\sin\phi_0$ is reasonably small.
So the old requirement $\xi_{k+1}\leq\xi<1$ is
stringent and not necessary for the quadratic asymptotic convergence.
Similar comments can be made on (\ref{suffcond2}) in Theorem~\ref{conds}
as well. We should point out that the above bound
for $\xi$ is conservative, and numerical experiments will demonstrate
that the inexact RQI with Lanczos works well and achieves
the quadratic convergence when $\xi$ exceeds bound (\ref{xichoice}).

{\bf Remark 2}. Conditions (\ref{linconv}) and (\ref{linconv2}) for
linear convergence show that $\xi_{k+1}$ can be as big as
$O(\frac{1}{\sin\phi_k})$ and $O(\frac{1}{\|r_k\|})$ as outer iterations
proceed. More precisely, (\ref{linconv}) indicates that
the inexact RQI with Lanczos still
converges linearly even if $\xi_{k+1}$ is as big as
$O(\frac{1}{\sin\phi_k})$ with the order constant
$\frac{\gamma}{2\beta+\gamma}$ smaller than one.

As done in \cite{jia09,SimonciniRQI}, we now estimate
$\|w_{k+1}\|$ in (\ref{EqIRQILinearEquation1}) obtained by Lanczos.
Note that the exact solution of
$(A-\theta_k I)w=u_k$ is $w_{k+1}=(A-\theta_k I)^{-1}u_k$.
Therefore, setting $\xi_{k+1}=0$ in (\ref{EqIRQIwk}), we have
\begin{eqnarray*}
\|w_{k+1}\|&=&\frac{\cos\phi_k}{|\theta_k-\lambda_1|}+O(\sin\phi_k)\\
&\approx&\frac{1}{|\theta_k-\lambda_1|}=\|(A-\theta_kI)^{-1}\|\\
&=&O\left(\frac{1}{\sin^2\phi_k}\right)=O\left(\frac{1}{\|r_k\|^2}\right),
\end{eqnarray*}
the last equality being from (\ref{parlett}).
From (\ref{EqIRQIwk}) and (\ref{parlett}), we also see that
these estimates hold for $\xi_{k+1}=O(\sin\phi_k)=O(\|r_k\|)$. So $\|w_{k+1}\|$
is also $O(\frac{1}{\|r_k\|^2})$ when the inexact RQI with Lanczos converges
cubically. Next we derive quantitative estimates on $\|w_{k+1}\|$ under more general
conditions, and by the estimates we establish a new quadratic convergence result.

\begin{theorem}\label{simplelanczos}
Let $\xi$ be a constant such that $\xi_{k+1}\leq\xi$ and satisfy
$\xi\sin\phi_k\leq\alpha<1$ with $\alpha$ a given constant not near one.
Then we have
\begin{eqnarray}
\|w_{k+1}\|&\geq&\frac{(1-\alpha)|\lambda_2-\lambda_1|}{4\beta\|r_k\|^2},
\label{wk+1lanczos}\\
\|r_{k+1}\|&\leq&\frac{\sqrt{1+\xi^2}}{\|w_{k+1}\|},\label{innerouter}\\
\|r_{k+1}\|&\leq&\frac{4\beta\sqrt{1+\xi^2}}{|\lambda_2-\lambda_1|(1-\alpha)}
\|r_k\|^2,\label{resrelanczos}
\end{eqnarray}
where {\rm (\ref{wk+1lanczos})} and {\rm (\ref{resrelanczos})} hold asymptotically
and {\rm (\ref{innerouter})}
holds exactly.
\end{theorem}

\begin{proof}
As $\xi_{k+1}d_{k+1}\perp u_k$ and $(A-\theta_k I)w_{k+1}=u_k+\xi_{k+1}d_{k+1}$, we
have
$$
\|(A-\theta_kI)w_{k+1}\|^2=\|u_k\|^2+\|\xi_{k+1}d_{k+1}\|^2=1+\xi_{k+1}^2.
$$
So we get from  $u_{k+1}=w_{k+1}/\|w_{k+1}\|$ that
\begin{equation}
\|(A-\theta_k
I)u_{k+1}\|=\frac{\sqrt{1+\xi_{k+1}^2}}{\|w_{k+1}\|}\leq
\frac{\sqrt{1+\xi^2}}{\|w_{k+1}\|}.\label{lanczosrela}
\end{equation}
By the optimality of Rayleigh quotient we obtain
\begin{equation}
\|r_{k+1}\|=\|(A-\theta_{k+1}I)u_{k+1}\|\leq\|(A-\theta_kI)u_{k+1}\|
=\frac{\sqrt{1+\xi_{k+1}^2}}{\|w_{k+1}\|}\leq\frac{\sqrt{1+\xi^2}}{\|w_{k+1}\|},
\label{optimal}
\end{equation}
which shows (\ref{innerouter}). It is easy to verify (cf.
\cite[p. 77]{ParlettSEP}) that
\begin{equation}
|\lambda_2-\lambda_1|\sin^2\phi_k\leq |\lambda_1 - \theta_k| \leq
(\lambda_{\max}-\lambda_{\min})\sin^2\phi_k. \label{error1}
\end{equation}

By using (\ref{EqIRQIwk}),
(\ref{denolanczos}), (\ref{error1}) and (\ref{parlett}) in turn,
we obtain
\begin{eqnarray*}
\|w_{k+1}\|&\geq&\frac{|\cos\phi_k+\xi_{k+1}\cos\psi_k|}{|\theta_k-\lambda_1|} \nonumber\\
&\geq& \frac{|1-\xi_{k+1}\sin\phi_k|}{|\theta_k-\lambda_1|} \label{lower}\\
&\geq&\frac{1-\alpha}{|\theta_k-\lambda_1|} \nonumber \\
&\geq&\frac{1-\alpha}{(\lambda_{\max}-\lambda_{\min})\sin^2\phi_k} \nonumber\\
&\geq&\frac{(1-\alpha)|\lambda_2-\lambda_1|}{4\beta\|r_k\|^2} \nonumber \\
&\geq&\frac{(1-\alpha)|\lambda_2-\lambda_1|}{4\beta\|r_k\|^2},\nonumber
\end{eqnarray*}
which proves (\ref{wk+1lanczos}). Substituting
(\ref{wk+1lanczos}) into (\ref{innerouter})
establishes (\ref{resrelanczos}).
\end{proof}

(\ref{resrelanczos}) indicates that the
inexact RQI with Lanczos converges quadratically if $\alpha$
is not near one. By combining this theorem with
Theorems~\ref{convlanczos}--\ref{conds},
(\ref{wk+1lanczos}) shows that $\|w_{k+1}\|$ is
always no less than $O(\frac{1}{\|r_k\|^2})$ when the inexact RQI
with Lanczos converges quadratically provided
that $\alpha$ is not near one. Noting that the inexact RQI with Lanczos
for $\xi_{k+1}=O(\sin\phi_k)=O(\|r_k\|)$ converges cubically
and $\|w_{k+1}\|$ is also $O(\frac{1}{\|r_k\|^2})$ (cf. the comments
before Theorem~\ref{simplelanczos}),
this illustrates that the size of $\|w_{k+1}\|$ itself
cannot reveal cubic and quadratic convergence rates of the inexact
RQI with Lanczos.  Furthermore, we cannot recover the cubic convergence of
the exact RQI and the inexact RQI when $\xi_{k+1}=\xi=0$ and
$\xi_{k+1}=\xi=O(\sin\phi_k)=O(\|r_k\|)$, respectively. So (\ref{resrelanczos}) is weaker
than Theorems~\ref{convlanczos}--\ref{conds}. This is because
$\|r_{k+1}\|\leq\|(A-\theta_kI)u_{k+1}\|$ is not sharp in the proof.

So far, all the convergence results are local, that is, they care
how the exact and inexact RQI behaves only from the current outer
iteration to the next one, assuming that current $(\theta_k,u_k)$ is
already a reasonably good approximation to $(\lambda,x)$. As is well
known, one of the important properties of the exact RQI is its {\em
global} residual monotonic decreasing property, i.e.,
$\|r_{k+1}\|\leq\|r_k\|$, for any (poor) starting vector $u_0$; see
Theorem 4.8.1 of \cite[p. 79]{ParlettSEP}. We now present a global
property to the inexact RQI with Lanczos.

\begin{theorem}\label{globallanczos}
For the inexact RQI with Lanczos starting with any starting vector
$u_0$, we have
\begin{equation}
\|r_{k+1}\|\leq\sqrt{1+\xi_{k+1}^2}\|r_k\|,\
k=0,1,\ldots.\label{monotonic1}
\end{equation}
\end{theorem}

\begin{proof}
From (\ref{EqIRQILinearEquation1}), we have
$$
w_{k+1}=(A-\theta_kI)^{-1}(u_k+\xi_{k+1}d_{k+1}).
$$
Again, note that for the Lanczos method its residual $\xi_{k+1}d_{k+1}$ satisfies
$\xi_{k+1}d_{k+1}^*u_k=0$. Then from (\ref{optimal}) and the Cauchy--Schwarz
inequality, we get
\begin{eqnarray*}
\frac{\|r_{k+1}\|}{\|r_k\|}&\leq&\frac{\|(A-\theta_kI)u_{k+1}\|}{\|r_k\|}
=\frac{\sqrt{1+\xi_{k+1}^2}}{\|w_{k+1}\|\|r_k\|}\\
&=&\frac{\sqrt{1+\xi_{k+1}^2}}{\|(A-\theta_kI)^{-1}(u_k+\xi_{k+1}d_{k+1})\|
\|(A-\theta_kI)u_k\|}\\
&\leq&\frac{\sqrt{1+\xi_{k+1}^2}}{|(u_k+\xi_{k+1}d_{k+1})^*(A-\theta_kI)^{-1}
(A-\theta_kI)u_k|}\\
&=&\frac{\sqrt{1+\xi_{k+1}^2}}{|1+\xi_{k+1}d_{k+1}^*u_k|}=\sqrt{1+\xi_{k+1}^2},
\end{eqnarray*}
which proves (\ref{monotonic1}).
\end{proof}

This theorem shows that, unlike the exact RQI, $\|r_k\|$
obtained by the inexact RQI with Lanczos is not monotonic
decreasing in the global sense for an arbitrary starting vector
$u_0$. This is similar to the inexact RQI with MINRES,
where Simoncini and Eld$\acute{e}$n~\cite{SimonciniRQI} have
derived a similar relationship between $\|r_{k+1}\|$ and $\|r_k\|$,
showing that residuals obtained by the inexact RQI with MINRES lose the
monotonic decreasing property that the exact RQI possesses; see Theorem 5.4
of \cite{SimonciniRQI}. Therefore,
as far as global convergence is concerned, the inexact RQI has a very
essential difference from the exact RQI, the former cannot guarantee its
convergence while the latter almost always converges for an arbitrary
starting vector \cite{ParlettSEP}. For the inexact RQI, we can only
expect its local convergence starting with a reasonably good starting vector.
We have made numerical experiments on some matrices for some
starting $u_0$'s generated randomly and found that it is indeed
the case for the inexact RQI with Lanczos or MINRES.

\section{Convergence of the inexact RQI with a tuned preconditioned
Lanczos}\label{precondit}

We have found that for a given $\xi$ satisfying our convergence conditions,
we may still need many inner iteration steps at each outer
iteration. This is especially true for difficult problems,
i.e., big $\beta$'s, or for computing
an interior eigenvalue $\lambda_1$ since it leads to a highly Hermitian
indefinite matrix $(A-\theta_k I)$ at each outer iteration. So,
in order to improve the overall performance, preconditioning
is generally necessary to speed up the Lanczos method. Some preconditioning
techniques have been proposed in e.g., \cite{mgs06,SimonciniRQI}. In
the unpreconditioned case, the right-hand side $u_k$ of
(\ref{EqERQILinearEquation1}) is rich in the direction of the
desired $x_1$. We can benefit much from this property when solving
the linear system. Actually, if the right-hand side is an eigenvector of the
coefficient matrix, Krylov subspace type methods will find the exact
solution in one step. However, a
usual preconditioner loses this important property, so that inner
iteration steps may not be reduced
\cite{mgs06,freitagspence,freitag08b}. A preconditioner with tuning
is necessary to recover this property and meanwhile attempts to improve the
conditioning of the preconditioned system, so that considerable
improvement over a usual preconditioner is possible
\cite{freitagspence,freitag08b,xueelman}. In what follows we show how to extend
our previous theory to the inexact RQI with a tuned preconditioned Lanczos.

Let $Q=LL^*$ be a Cholesky factorization of some Hermitian positive
definite matrix which is an approximation to $A-\theta_k I$ in some
sense \cite{mgs06,freitag08b,xueelman}. A tuned preconditioner
${\cal Q}={\cal L}{\cal L}^*$ can be constructed by adding a rank-1
or rank-2 modification to $Q$, so that
\begin{equation}
{\cal Q}u_k=Au_k;\label{tune}
\end{equation}
see \cite{freitagspence,freitag08b,xueelman} for details. Using
the tuned preconditioner ${\cal Q}$, the shifted inner linear system
(\ref{EqERQILinearEquation1})
is equivalently transformed to the preconditioned one
\begin{equation}
B\hat{w}={\cal L}^{-1}(A-\theta_k I){\cal L}^{-*}\hat{w}={\cal
L}^{-1}u_k \label{precond}
\end{equation}
with the original $w={\cal L}^{-*}\hat{w}$. Once the Lanczos method
is used to solve it, we are led to the inexact RQI with a tuned
preconditioned Lanczos. A power of the tuned preconditioner ${\cal
Q}$  is that the right-hand side ${\cal L}^{-1}u_k$ is rich in the
eigenvector of $B$ associated with its smallest eigenvalue and has
the same quality as $u_k$ as an approximation to the eigenvector $x_1$ of
$A$, while for the usual preconditioner $Q$ the right-hand side
$L^{-1}u_k$ does not possess this property.

Take the zero vector as an initial guess to the solution of
(\ref{precond}) and let $\hat{w}_{k+1}$ be the approximate solution
obtained by the $m$-step Lanczos method applied to it. Then we have
\begin{equation}
{\cal L}^{-1}(A-\theta_k I){\cal L}^{-*}\hat{w}_{k+1}={\cal
L}^{-1}u_k+\hat{\xi}_k\hat{d}_{k+1}, \label{preresidual}
\end{equation}
where $\hat{w}_{k+1}\in {\cal K}_m(B,{\cal L}^{-1}u_k)$,
$\hat{\xi}_{k+1}\hat{d}_{k+1}$ with $\|\hat{d}_{k+1}\|=1$ is the residual and
$\hat{d}_{k+1}$ is the residual direction vector. Keep in mind that
$w_{k+1}={\cal L}^{-*}\hat{w}_{k+1}$. We then get
\begin{equation}
(A-\theta_k I)w_{k+1}=u_k+\hat{\xi}_{k+1}{\cal
L}\hat{d}_{k+1}=u_k+\hat{\xi}_{k+1}\|{\cal L}\hat{d}_{k+1}\|\frac{{\cal
L}\hat{d}_{k+1}}{\|{\cal L}\hat{d}_{k+1}\|}.\label{orig}
\end{equation}
So $\xi_{k+1}$ and $d_{k+1}$ in (\ref{EqIRQILinearEquation1}) are
$\hat{\xi}_{k+1}\|{\cal L}\hat{d}_{k+1}\|$ and $\frac{{\cal
L}\hat{d}_{k+1}}{\|{\cal L}\hat{d}_{k+1}\|}$, respectively. Hence our
general Theorems~\ref{ThmIRQIQuadraticConvergence}--\ref{resbound}
apply and are not repeated here.

An extension of Theorem~\ref{lanczos} to the preconditioned case is
nontrivial and needs more work.
Let $(\mu_i,y_i),\ i=1,2,\ldots,n$ be the eigenpairs of $B$ with
$$
\mid\mu_1\mid<\mid\mu_2\mid\leq\cdots\leq\mid\mu_n\mid.
$$
Define $\hat{u}_k={\cal L}^{-1}u_k/\|{\cal L}^{-1}u_k\|$. Similar to
{\rm (\ref{EqIRQIDecompositoinOfu_k})} and {\rm
(\ref{EqIRQIDecompositoinOfd_{k+1}}), let
\begin{eqnarray}
\hat{u}_k&=&y_1\cos\hat{\phi}_k+\hat{e}_k\sin\hat{\phi}_k,
\ \hat{e}_k\perp y_1, \ \|\hat{e}_k\|=1,
\label{predec1}\\
\hat{d}_{k+1}&=&y_1\cos\hat{\psi}_k+\hat{f}_k\sin\hat{\psi}_k,
\ \hat{f}_k\perp y_1, \ \|\hat{f}_k\|=1 \label{predec2}
\end{eqnarray}
be the orthogonal direct sum decompositions. Then it is
known \cite{freitag08b} that
\begin{eqnarray}
|\mu_1|&=&O(\sin\phi_k),\label{mu1}\\
\sin\hat{\phi}_k&\leq& c_1\sin\phi_k \label{hatphi}
\end{eqnarray}
with $c_1$ a constant.




Similar to Theorem~\ref{lanczos}, we can derive the following results.

\begin{theorem}\label{planczos3}
It holds that
\begin{eqnarray}
|\cos\psi_k|&=&O(\sin\phi_k),\label{planczos4}\\
\sin\psi_k&=&1-O(\sin^2\phi_k).\label{planczos5}
\end{eqnarray}
\end{theorem}

\begin{proof}

For the $m$-step Lanczos method for
(\ref{preresidual}), we have $\hat{\xi}_{k+1}\hat{d}_{k+1}\perp {\cal
K}_m(B,{\cal L}^{-1}u_k)$. Particularly, it holds that
$$
\hat{d}_{k+1}^*{\cal L}^{-1}u_k=\hat{d}_{k+1}^*{\cal L}^{-1}u_k=0,
$$
from which and $\hat{u}_k={\cal L}^{-1}u_k/\|{\cal L}^{-1}u_k\|$ it follows that
$\hat{d}_{k+1}^*\hat{u}_k=0$.
Therefore, from $d_{k+1}=\frac{{\cal L}\hat{d}_{k+1}}{\|{\cal
L}\hat{d}_{k+1}\|}$ we have
$$
0=\hat{d}_{k+1}^*\hat{u}_k=\hat{d}_{k+1}^*\underbrace{{\cal L}^*{\cal L}^{-*}}_I\hat{u}_k
=\|{\cal L}\hat{d}_{k+1}\|\|{\cal L}^{-1}u_k\|d_{k+1}^*{\cal L}^{-*}\hat{u}_k=0,
$$
i.e.,
\begin{equation}
d_{k+1}^*{\cal L}^{-*}\hat{u}_k=0.\label{pdec}
\end{equation}

By definition, we have
$$
{\cal L}^{-1}(A-\theta_k I){\cal L}^{-*}y_1=\mu_1y_1,
$$
from which it follows that
$$
(A-\theta_kI)\tilde{u}_k=\mu_1\frac{{\cal L}y_1}{\|{\cal L}^{-*}y_1\|}
$$
with $\tilde{u}_k={\cal L}^{-*}y_1/\|{\cal L}^{-*}y_1\|$. Therefore,
by standard perturbation theory and (\ref{mu1}), we get
\begin{equation}
\sin\angle(\tilde{u}_k,x)=O(\mu_1\frac{\|{\cal L}y_1\|}
{\|{\cal L}^{-*}y_1\|})=O(|\mu_1|)=O(\sin\phi_k). \label{tildeuk}
\end{equation}
On the other hand, from (\ref{hatphi}), we have
$$
\sin\hat{\phi}_k=\sin\angle(\hat{u}_k,y_1)= \sin\angle({\cal
L}^{-1}u_k,y_1)=O(\sin\phi_k).
$$
Therefore, we can write
$$
\hat{u}_k=y_1+O(\sin\phi_k),
$$
which leads to
$$
{\cal L}^{-*}\hat{u}_k={\cal L}^{-*}y_1+O(\sin\phi_k).
$$
Thus, we have
\begin{equation}
\sin\angle({\cal L}^{-*}\hat{u}_k,{\cal L}^{-*}y_1)
=\sin\angle({\cal L}^{-*}\hat{u}_k,\tilde{u}_k)=O(\sin\phi_k).\label{hatuk}
\end{equation}

Since
$$
\angle({\cal L}^{-*}\hat{u}_k,x)\leq\angle({\cal L}^{-*}\hat{u}_k,
\tilde{u}_k)+\angle(\tilde{u}_k,x),
$$
combining (\ref{tildeuk}) and (\ref{hatuk}), we get
\begin{equation}
\sin\angle({\cal L}^{-*}\hat{u}_k,x)\leq\sin\angle({\cal L}^{-*}\hat{u}_k,
\tilde{u}_k)+\sin\angle(\tilde{u}_k,x)=O(\sin\phi_k).\label{lhatuk}
\end{equation}
Recall that $d_{k+1}=x\cos\psi_k+e_k\sin\psi_k$ and substituting it and
the orthogonal direct sum decomposition
$$
{\cal L}^{-*}\hat{u}_k=\|{\cal L}^{-*}\hat{u}_k\|
(x\cos\angle({\cal L}^{-*}\hat{u}_k,x)+g_k\sin\angle({\cal L}^{-*}\hat{u}_k,x))
$$
with $g_k\perp x$ into (\ref{pdec}). Then following the proof of Theorem~\ref{lanczos},
we can get
\begin{eqnarray*}
|\cos\psi_k|&\leq&|\tan\angle({\cal L}^{-*}\hat{u}_k,x)|,\\
\sin\psi_k&=&1-O(\sin^2\angle({\cal L}^{-*}\hat{u}_k,x)).
\end{eqnarray*}
Combining them with (\ref{lhatuk}) yields (\ref{planczos4}) and (\ref{planczos5}).

\end{proof}

Using this theorem and writing (\ref{planczos4}) as
$|\cos\psi_k|\leq c_2\sin\phi_k$ with $c_2$ a constant, it is direct
to extend Theorems~\ref{convlanczos}--\ref{conds} in the
unpreconditioned Lanczos case to the tuned preconditioned Lanczos
case. We have done preliminary numerical experiments and confirmed
the theory. Our concerns in this paper are only the convergence
theory of the inexact RQI with the unpreconditioned and tuned
preconditioned Lanczos, and the pursue of effective tuned
preconditioners is beyond the scope of the current paper. We will
only report numerical results on the inexact RQI with the
unpreconditioned Lanczos.

\section{Numerical experiments}\label{testlanczos}

Our numerical experiments were performed on an Intel
(R) Core (TM)2 Quad CPU Q9400 $2.66$GHz with main memory 2 GB using
Matlab 7.8.0 with the machine precision $\epsilon=2.22\times
10^{-16}$ under the Microsoft Windows XP operating system.

We report the numerical results by the inexact RQI with the
unpreconditioned Lanczos for computing the smallest eigenpairs
of four symmetric (Hermitian)
matrices: BCSPWR08 of order 1624, CAN1054 of order 1054, DWT2680 of
order 3025 and LSHP3466 of order 3466 \cite{duff}. Recall the definition
(\ref{beta}) of $\beta$. Note that the
bigger the factor $\beta$ is, the worse conditioned $x_1$ is. Meanwhile,
for $\beta$ big, Theorem~\ref{ThmIRQIQuadraticConvergence} and
Theorem~\ref{convlanczos} show that although RQI and the inexact RQI
can still converge cubically and quadratically, they may converge
more slowly and needs more outer iterations as the
factors $2\beta$ and $\frac{2\beta\xi}{1-\alpha}$ in
(\ref{bound1}) and (\ref{lanczosquadra2}) are big. As a reference, we
use the Matlab function {\sf eig.m} to compute $\beta$. We find that
DWT2680 and LSHP3466 are considerably more difficult than the
other two. We only report the results on the computation of
the smallest eigenpair.

Theorems~\ref{convlanczos}--\ref{conds} tells us
that the cubic asymptotic convergence of the inexact RQI with
Lanczos is achieved for $\xi_{k+1}=O(\sin\phi_k)=O(\|r_k\|)$
when updating $(\theta_k,u_k)$
to get $(\theta_{k+1},u_{k+1})$, in the experiments we take
\begin{equation}
\xi_{k+1}\leq \frac{\|r_k\|}{\|A\|_1}.\label{decrease}
\end{equation}
Other stopping criteria have been taken, e.g.,
$\xi_{k+1}\leq\min\{\tau,\tau \|r_k\|\}$ with $\tau=0.1$ in \cite{mgs06}. They
are essentially the same as (\ref{decrease}) and differ only with the scaling
factor before $\|r_k\|$. But (\ref{decrease}) may be more general
as it takes the size of $A$ into account.

We construct the same initial $u_0$ for each matrix that is $x_1$ plus
a reasonably small perturbation generated randomly in a uniform
distribution, such that
$|\lambda_1-\theta_0|<\frac{|\lambda_2-\lambda_1|}{2}$. The algorithm
stops whenever $\|r_k\|=\|(A-\theta_k I)u_k\|\leq\|A\|_1tol$, where
$tol=10^{-14}$ unless stated otherwise. In the experiments, we
use the Matlab function {\sf symmlq.m} to solve the inner linear
systems when $\xi_k<1$. We should notice that for $\xi_k\geq
1$ the Matlab function {\sf symmlq.m} cannot be applied. Since the
Lanczos method behaves irregularly and may nearly break down or
break down for indefinite linear systems, that is, $T_m$ in
(\ref{lanczos}) is ill conditioned and can be nearly singular and
even numerically singular, it may produce bad approximate solutions
with large norms and large residual norms $\xi_k$'s for some steps
$m$. As far as solving the linear systems is concerned, such
approximate solutions have no accuracy and no practical value. In
{\sf symmlq.m}, if such a bad approximate solution emerges, it
always outputs the approximate solution as zero and the residual
norm $\xi_k=1$ simply, telling us nothing! However, we have seen
that in the inexact RQI with Lanczos, $\xi_k\geq 1$ is allowed. So
for our purpose, we have worked out a Lanczos code that uses the
Gram--Schmidt with iterative refinement \cite{stewart} to generate a
numerically orthonormal basis of the Krylov subspace ${\cal
K}_m(A,u_k)$ and delivers 'correct' results that the Lanczos method
should produce. We point out that our Lanczos code is not optimized
but numerically stable.

We report the results obtained by the inexact RQI with Lanczos
for choosing $\xi_{k+1}$ as in
(\ref{decrease}) and fixed $\xi=0.1,1,5$. Based on
our theory, the method should asymptotically converge quadratically for the
$\xi=1,5$ and use almost the same outer iterations as those for
$\xi=0.1$. Therefore, the total computational cost may be reduced
considerably. For the inexact RQI with Lanczos, the total inner
iteration steps ``$iters$'', i.e., the total matrix-vector products in inner
iterations, is a good and reasonable
measure of overall performance of the method, as commonly adopted in many cited
papers, e.g., \cite{SimonciniRQI} and those of Spence and his coworkers.

Tables~\ref{BCSPWR08lanczos}--\ref{lshplanczos} list the computed
results, where $iters$ denotes the number of total inner iteration
steps, $iter^{(k)}$ the number of inner iteration steps at the $k$-th
outer iteration and the "-" denotes the stagnation of {\sf
symmlq.m} at the $iter^{(k)}$-th step. We
comment that in {\sf symmlq.m} the output $iter^{(k)}=m-1$, where
$m$ is the steps of the Lanczos process.

\begin{table}[htp]
\begin{center}
\begin{tabular}{|c|c|c|c|c|c|c|}\hline
$\xi_k\leq\xi$&$k$&$\|r_k\|$&$\sin\phi_k$&$\xi_k$&$iter^{(k)}$
&$iters$\\\hline
0 (RQI)&1& 0.0124&0.0036& & &  \\
&  2&$3.1e-8$&$8.5e-8$& & & \\
& 3&$2.0e-15$&$6.6e-15$ & & &\\\hline
$\frac{\|r_{k-1}\|}{\|A\|_1}$&1 &0.0071&0.0029&0.0367&7&1003\\
&2&$3.7e-8$&$1.2e-7$&$4.1e-4$&40&\\
&3&$3.0e-15$&$3.0e-15$&-&956&\\\hline
0.1&1&0.0090&0.0045&0.0950&5&300\\
&2&$2.1e-5$&$3.7e-5$&0.0847&24&\\
&3&$3.3e-13$&$3.2e-13$&0.0754&44&\\\hline
1&1&0.0462&0.0165&0.7161&3&87\\
&2&$5.1e-4$&$8.9e-4$&0.8505&11&\\
&3&$2.3e-7$&$2.9e-7$&0.9829&26&\\
&4&$1.1e-14$&$1.3e-14$&-&47&\\\hline
5&1&0.1259&0.0332&2.0694&2&87\\
&2&0.0107&0.0064&3.0112&4&\\
&3&$1.9e-4$&$2.5e-4$&4.2267&13&\\
&4&$1.2e-7$&$1.6e-7$&4.8117&25&\\
&5&$5.7e-14$&$3.5e-14$&4.8558&43&\\\hline
\end{tabular}
\begin{tabular}{|c|c|c|}\hline
$m$ & outer iterations& $iters$\\\hline
5 &110 &550\\
10 &21&210\\
15& 10 &150\\
20 &7 &140\\
30 & 5 &150\\\hline
\end{tabular}
\caption{BCSPWR08, $\beta=40.19,\
\sin\phi_0=0.1020$.}\label{BCSPWR08lanczos}
\end{center}
\end{table}

\begin{table}[ht]
\begin{center}
\begin{tabular}{|c|c|c|c|c|c|c|}\hline
$\xi_k\leq\xi$&$k$&$\|r_k\|$&$\sin\phi_k$&$\xi_k$&
$iter^{(k)}$&$iters$\\\hline
0 (RQI)&1&0.0269&0.0110& & &  \\
& 2&$1.1e-7$&$2.1e-7$ & & & \\
& 3&$3.2e-15$&$5.0e-15$ & &  &\\
&4&$2.5e-15$&$4.7e-15$ & &  &\\\hline
$\frac{\|r_{k-1}\|}{\|A\|_1}$&1 &0.0155&0.0038&0.0403&8&432\\
&2&$1.4e-7$&$3.1e-7$&$3.1e-4$&45&\\
&3&$1.4e-14$&$5.0e-15$&-&379&\\\hline
0.1&1&0.0181&0.0055&0.0796&6&207\\
&2&$3.8e-6$&$2.2e-5$&0.0679&30&\\
&3&$4.7e-13$&$2.0e-13$&0.0969&51&\\
&4&$6.0e-15$&$5.0e-15$&-&120&\\\hline
1&1&0.0624&0.0132&0.4376&4&107\\
&2&$5.7e-4$&$9.7e-4$&0.9980&14&\\
&3&$2.4e-7$&$1.3e-7$&0.9561&30&\\
&4&$1.7e-14$&$8.0e-14$&0.9952&53&\\\hline
5&1&0.2238&0.0454&1.7072&2&101\\
&2&0.0310&0.00946&3.7008&4&\\
&3&$8.3e-4$&0.0012&4.5335&11&\\
&4&$1.7e-6$&$1.1e-6$&3.6374&26&\\
&5&$4.0e-12$&$1.5e-12$&3.9095&45&\\
&6&$3.4e-15$&$2.0e-14$&4.5305&13&\\\hline
\end{tabular}
\begin{tabular}{|c|c|c|}\hline
$m$ & outer iterations& $iters$\\\hline
5 &110 &550\\
10 &21&210\\
15& 10 &150\\
20 &7 &140\\
30 & 5 &150\\\hline
\end{tabular}
\caption{CAN1054, $\beta=88.28,\
\sin\phi_0=0.1008$.}\label{can1054lanczos}
\end{center}
\end{table}

\begin{table}[ht]
\begin{center}
\begin{tabular}{|c|c|c|c|c|c|c|}\hline
$\xi_k\leq\xi$&$k$&$\|r_k\|$&$\sin\phi_k$&$\xi_k$&$iter^{(k)}$
&$iters$\\\hline
0 (RQI)&1&0.0144&0.1188& & &  \\
&  2&$1.6e-4$&0.0018& &  & \\
& 3&$5.1e-10$&$2.6e-8$ & & &\\
&4&$1.0e-15$&$6.7e-13$ & & & \\\hline
$\frac{\|r_{k-1}\|}{\|A\|_1}$&1 &0.0143&0.0121&0.1171&5&1512\\
&2&$4.8e-5$&$5.1e-4$&$0.0019$&167&\\
&3&$4.2e-11$&$9.9e-10$&-&424&\\
&4&$6.1e-15$&$6.6e-13$&-&916&\\\hline
0.1&1&0.0123&0.0104&0.0830&6&955\\
&2&$5.1e-5$&$1.1e-4$&0.0972&93&\\
&3&$1.3e-11$&$1.5e-10$&0.0986&267&\\
&4&$5.0e-15$&$6.6e-13$&0.1109(*)&589&\\ \hline
1&1&0.0419&0.0265&0.4842&3&402\\
&2&$6.7e-4$&$0.0041$&0.9424&15&\\
&3&$8.6e-7$&$1.2e-5$&0.9667&117&\\
&4&$1.9e-12$&$2.2e-11$&0.9276&267&\\ \hline
5&1&0.0851&0.0387&1.0223&2&535\\
&2&0.0092&0.0125&3.8250&5&\\
&3&$2.8e-4$&0.0027&4.8280&19&\\
&4&$1.1e-6$&$1.3e-5$&4.7719&108&\\
&5&$1.3e-11$&$2.6e-10$&4.8329&238&\\
&6&$7.0e-15$&$5.8e-13$&4.9521&163&\\\hline
\end{tabular}
\begin{tabular}{|c|c|c|}\hline
$m$ & outer iterations& $iters$\\\hline
10 &268 &2680\\
20 &55&1100\\
30& 30 &900\\
40 &17 &680\\
50 & 11 &550\\
60 &11 & 660\\\hline
\end{tabular}
\caption{DWT2680, $tol=10^{-12}$, $\beta=2295.6,\
\sin\phi_0=0.1095$.}\label{dwtlanczos}
\end{center}
\end{table}

\begin{table}[ht]
\begin{center}
\begin{tabular}{|c|c|c|c|c|c|c|}\hline
$\xi_k\leq\xi$&$k$&$\|r_k\|$&$\sin\phi_k$&$\xi_k$&
$iter^{(k)}$&$iters$\\\hline
0 (RQI)&1&0.0149&0.1716& & &   \\
&  2&$4.0e-4$&0.0056& & &  \\
& 3&$1.2e-8$&$8.9e-7$ & & & \\
&4&$2.0e-15$&$4.0e-13$ & & &\\\hline
$\frac{\|r_{k-1}\|}{\|A\|_1}$&1 &0.0102&0.0097&0.0874&6&1717\\
&2&$3.9e-5$&$1.2e-4$&$0.0014$&201&\\
&3&$1.2e-9$&$4.5e-8$&-&497&\\
&4&$5.2e-15$&$6.1e-13$&-&1013&\\\hline
0.1&1&0.0102&0.0098&0.0874&6&651\\
&2&$4.2e-4$&$5.4e-4$&0.0990&102&\\
&3&$4.5e-8$&$2.8e-7$&0.0948&256&\\
&4&$5.3e-13$&$5.6e-11$&0.0965&287&\\ \hline
1&1&0.0408&0.0251&0.5335&3&424\\
&2&$6.4e-4$&0.0036&0.9685&15&\\
&3&$7.6e-7$&$1.9e-5$&0.9770&123&\\
&4&$2.7e-12$&$2.9e-11$&0.9907&283&\\\hline
5&1&0.0779&0.0370&1.0528&2&444\\
&2&0.0088&0.0128&4.2424&5&\\
&3&$2.4e-4$&0.0021&4.4646&21&\\
&4&$7.0e-7$&$2.8e-5$&4.6415&121&\\
&5&$1.9e-11$&$2.4e-10$&4.6506&146&\\
&6&$1.0e-14$&$3.8e-13$&4.6805&149&\\\hline
\end{tabular}
\begin{tabular}{|c|c|c|}\hline
$m$ & outer iterations& $iters$\\\hline
10 &278 &2780\\
20 &54&1080\\
30& 32 &960\\
40 &19 &760\\
50 & 13 &650\\
60 &11 & 660\\\hline
\end{tabular} \caption{LSHP3466, $tol=10^{-12}$,
$\beta=2613.1,\ \sin\phi_0=0.1011$.}\label{lshplanczos}
\end{center}
\end{table}

Before commenting the experiments, we should remind that in finite
precision arithmetic $\|r_k\|/\|A\|_1$ can not decline further
whenever it reaches a moderate multiple of $\epsilon=2.2\times
10^{-16}$. Therefore, assuming that the algorithm stops at outer
iteration $k$, if $\sin\phi_{k-1}$ or $\|r_{k-1}\|\leq 10^{-6}$ or
$10^{-9}$, then the algorithm may not continue converging cubically
or quadratically at the final outer iteration $k$. Another point is
that when judging convergence rates, we must take the factor $\beta$
into account. The smaller it is, the more clearly cubic and
quadratic convergence exhibits, as indicated by
(\ref{lanczosquadra}) and (\ref{resquadra}); the bigger it is, the
less apparent cubic and quadratic convergence is. So we should
precisely base (\ref{lanczosquadra}) or (\ref{resquadra}) to judge
cubic and quadratic convergence of $\sin\phi_k$ or $\|r_k\|$. In
this sense, we see from
Tables~\ref{BCSPWR08lanczos}--\ref{lshplanczos} that the exact RQI
and the inexact RQI with Lanczos for decreasing $\xi_k$ converge
cubically and the method starts to converge quadratically for the
given fixed $\xi$'s after very few outer iterations. Strikingly, for
$\xi=1,5$, the method is much more efficient than that for
$\xi=0.1$ and for decreasing $\xi_k$; the method with $\xi\geq 1$ is
four times and twice as fast as that with $\xi=0.1$ for BCSPWR08 and
for CAN1054 and DWT2680, respectively. For LSHP3466, the gain is not
so great, but the method with $\xi\geq 1$ is still one and a half
times as fast as the method with $\xi=0.1$. With the fixed
$\xi=1,5$, it is always three to ten times as fast as that with
decreasing $\xi_k=O(\|r_{k-1}\|)$ for the four test matrices. It is seen
that for a bigger $\xi$ the method may need a little more outer
iterations but it does indeed converge quadratically and is in
agreement with quadratic convergence bound (\ref{lanczosquadra}).
Why the method with fixed bigger $\xi_k$'s converges a
little more slowly is due to the bigger convergence factor
$\frac{2\beta\xi}{1-\alpha}$ in
(\ref{lanczosquadra2}).

Note that the linear systems $(A-\theta_k)w=u_k$'s are Hermitian
indefinite and become increasingly worse conditioned and even
numerically singular as $\theta_k\rightarrow\lambda$ with increasing
$k$. So, more inner iteration steps are needed generally for a fixed
$\xi$ as $k$ increases. We find that for the difficult DTW2680 and
LSHP3466, many more inner iterations are used than those for BCSPWR08 and
CAN1054.

We have tested many fixed $\xi$'s ranging from 10 to 50 for each
matrix and found that the method converges quadratically. $\xi\geq 10$ does not
satisfy condition (\ref{suffcond1}) for
quadratic convergence but is a moderate multiple of $\frac{1}{\sin\phi_0}$.
The algorithm with these $\xi$ behaves almost the same as that with $\xi=1,5$
and uses a little more outer iterations and comparable total inner
iteration steps $iters$.

Unlike quadratic convergence where it is not necessary to
estimate $\xi$ accurately, for linear convergence,
we see that conditions (\ref{linconv}) and
(\ref{linconv2}) {\em heavily} depend on and are sensitive to the a-priori
$\beta$. So it appears impossible to design a practical criterion robustly and
reliably unless a good estimate on $\beta$ is available in advance.
Note that the convergence of the method allows $\xi_{k+1}$ to increase
up to $O(\frac{1}{\|r_k\|})$ as outer iterations proceed. Therefore,
we may implement the inexact RQI
with Lanczos for certain fixed inner iteration steps $m$'s. Doing so is
based on a not very stringent expectation that resulting $\xi_{k+1}$'s are not
too large and at least obey one of (\ref{linconv}) and (\ref{linconv2}).
Of course, if $T_m$ is too ill conditioned, it is possible for $\xi_{k+1}$
to be too large and exceed  bounds (\ref{linconv}) and (\ref{linconv2}).
We have tested several $m$'s for each test matrix, see
Tables~\ref{BCSPWR08lanczos}--\ref{lshplanczos} for results.
Figure~\ref{figfixed} displays convergence processes of the inexact
RQI with Lanczos for various fixed $m$'s.

We find that the method works very well and robustly. In terms of
$iters$, it is seen from the tables that the overall efficiency of
the inexact RQI with Lanczos for fixed small $m$'s is comparable to
that of the method with given fixed $\xi$'s, except $m=5$ for the
easy problems BCSPWR08 and CAN1054 and $m=10$ for
the relatively difficult problems DWT2680 and
LSHP3466. Furthermore, we observe that for BCSPWR08 and CAN1054 the inexact
RQI with Lanczos converges almost as fast as the exact RQI for
$m=30$, and for the difficult problems DWT2680 and LSHP3466 we need
to properly increase $m$ to achieve fast outer convergence. As
expected, it is not surprising from the figures that a five-step
and at most ten-step Lanczos method for the inner linear systems is
enough to ensure the convergence of the inexact RQI with Lanczos.
However, for the difficult BWT2680 and LSHP3466, the inexact RQI
with Lanczos fails to converge when $m=5$. The reason is that some
$\xi_{k+1}$'s are too big and violate the linear convergence conditions.
Finally, we should point out that although the Lanczos method with a smaller
$m$ usually produces a bigger $\xi_{k+1}$ and makes the inexact RQI
use more outer iterations, the total
inner iteration steps $iters$ may not increase.

\begin{figure}[ht]
\begin{center}
\includegraphics[scale=.48]{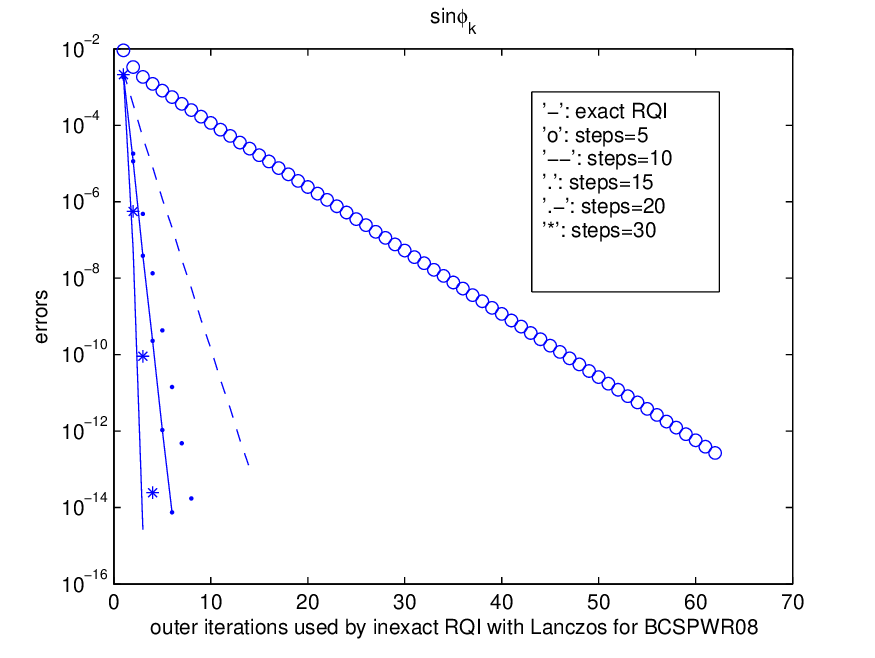}
\includegraphics[scale=.48]{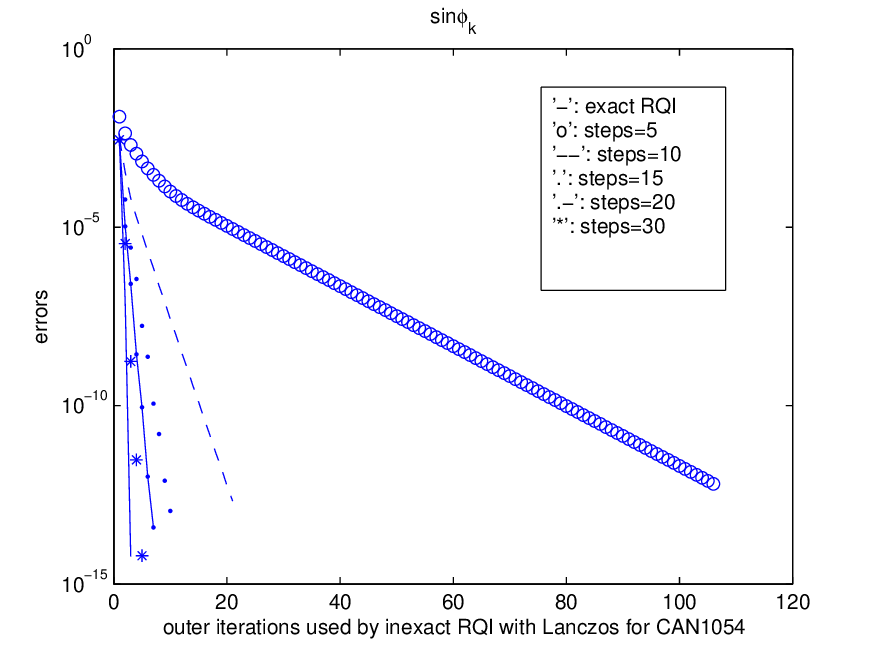}
\includegraphics[scale=.48]{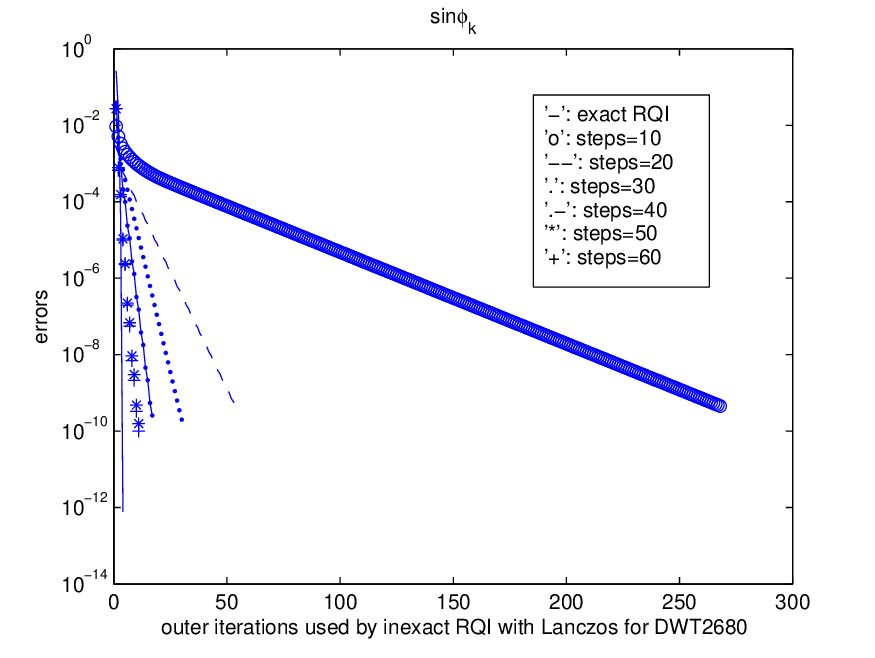}
\includegraphics[scale=.48]{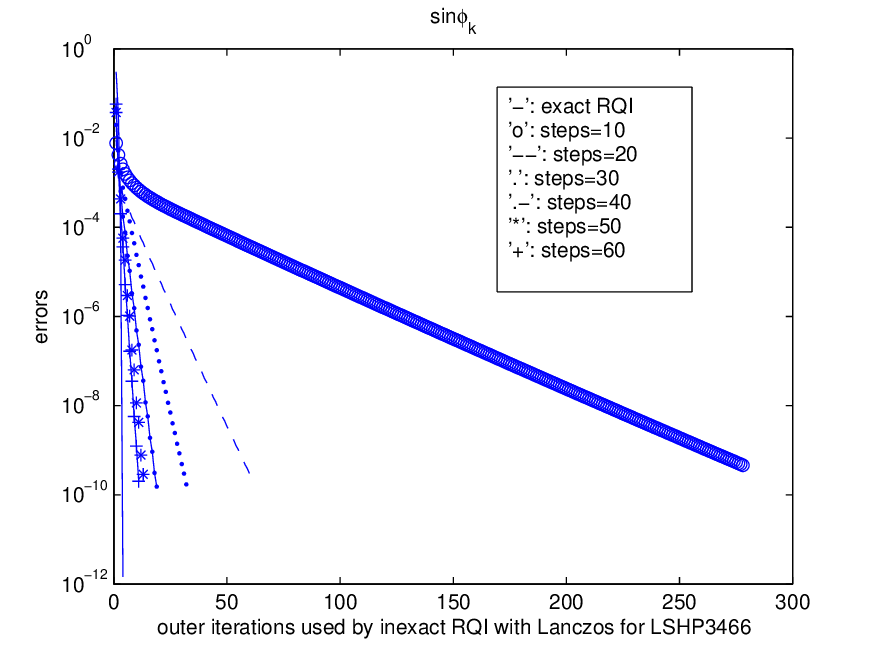}
\caption{The inexact RQI with Lanczos for varying inner iteration
steps} \label{figfixed}
\end{center}
\end{figure}

To see how big $\xi_{k+1}$'s may be as outer iterations converge, we
display the curves of $\xi_{k+1}$'s versus $\frac{1}{\sin\phi_k}$'s for
some fixed $m$'s in Figure~\ref{figbcspwrlanczos}. We remark that
$\xi_{k+1}$ in the figure should correspond to $\xi_k$ in our
context.

\begin{figure}[ht]
\begin{center}
\includegraphics[scale=.48]{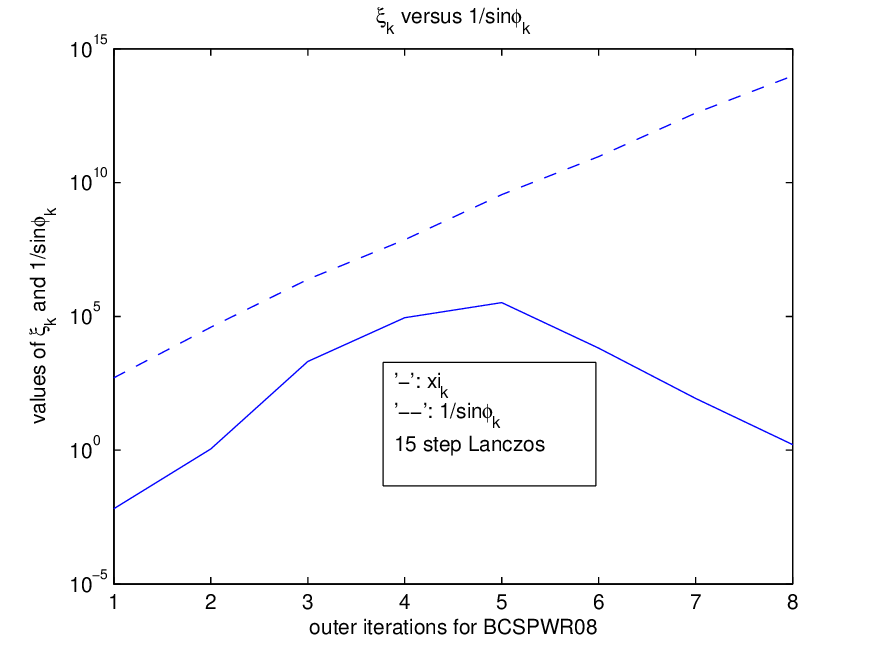}
\includegraphics[scale=.48]{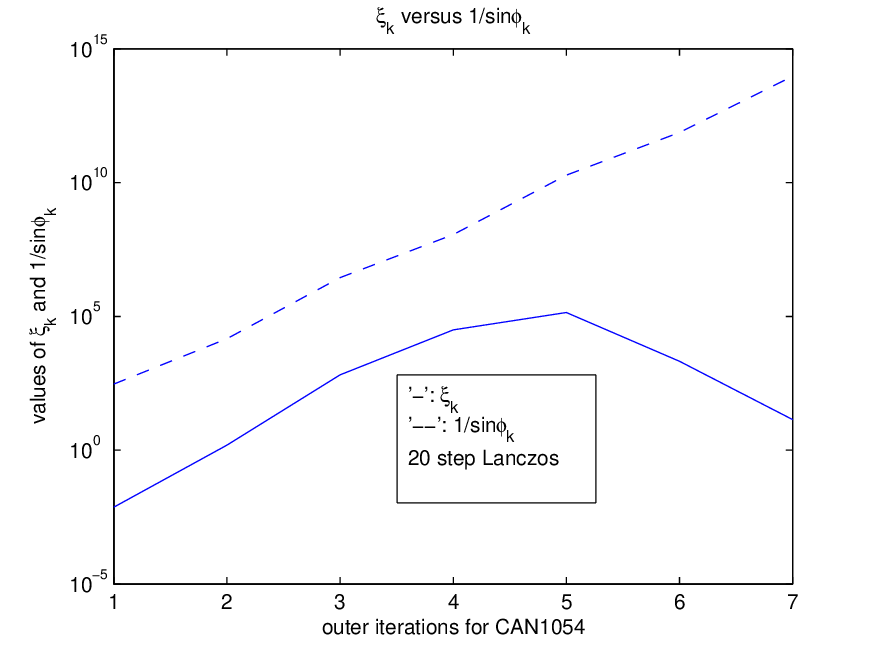}
\includegraphics[scale=.48]{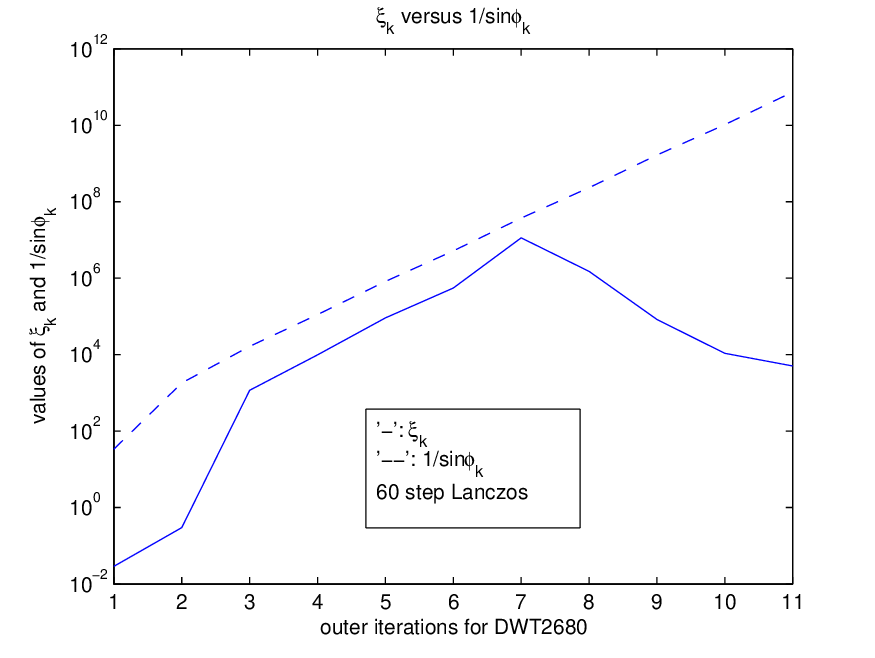}
\includegraphics[scale=.48]{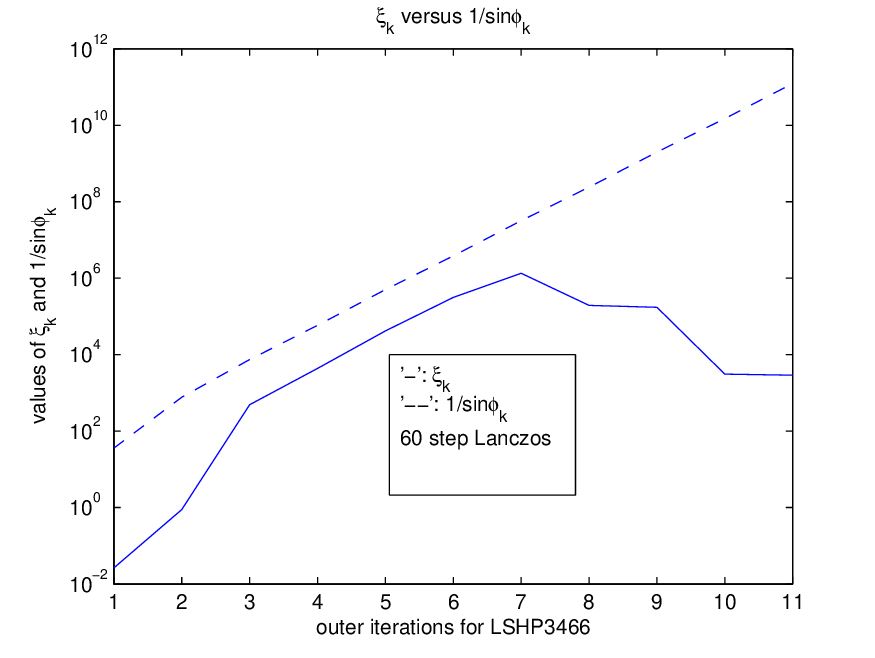}
\caption{$\xi_{k+1}$ versus $\frac{1}{\sin\phi_k}$. }
\label{figbcspwrlanczos}
\end{center}
\end{figure}

It is clear from Figure~\ref{figbcspwrlanczos} that the inexact RQI
with Lanczos works very well and $\sin\phi_k$ tends to zero smoothly
and quickly, that is, $\frac{1}{\sin\phi_k}$ tends to infinity
smoothly and quickly, though most of $\xi_{k+1}$'s are much bigger than
one and some of them are near $\frac{1}{\sin\phi_k}$ and can be as
big as $10^4\sim 10^{7}$! Furthermore, we observe that almost all
$\xi_{k+1}$'s are smaller than $\frac{1}{\sin\phi_k}$ and the only
exception is $\xi_6>\frac{1}{\sin\phi_6}$ for BWT2680.

We also computed some other eigenpairs of each test matrix. We
observed similar behavior and confirmed our theory. However, when an
interior eigenpair was required, we often needed much more $iters$.
This is because the shifted inner linear systems can be highly
indefinite (i.e., each shifted matrix has many positive and negative
eigenvalues) and may be hard to solve. If the smallest or largest
eigenpair is required, then $A-\theta_k I$ has only one negative or
positive eigenvalue, assuming that $A$ only has simple eigenvalues.
As a consequence, after the smallest Ritz value converges to the
smallest eigenvalue $\lambda-\theta_k$ of $A-\theta_k I$, the
Lanczos method will behave as if $A-\theta_k I$ is positive or
negative definite, so that it converges smoothly after the smallest
Ritz value converges. Comparing with the computation of the smallest
or largest eigenpair, we found that we must take a considerably
bigger fixed inner iteration steps $m$ to make the method converge
correctly when an interior eigenpair was desired. For a descriptive
analysis, see, e.g., \cite{freitagspence,freitag08b,paige,xueelman}.

\section{Conclusions}\label{conc}

We have considered the convergence of the inexact RQI with the
unpreconditioned and tuned preconditioned Lanczos methods and have
established a number of results. These results show how inner
tolerance affects accuracy of outer iterations and provide practical
criteria on how to best control inner tolerance to achieve the quadratic
asymptotic convergence of the inexact RQI. It is the first time to appear
surprisingly that the inexact RQI with Lanczos converge
quadratically provided $\xi_{k+1}\leq\xi$ with $\xi$ a constant being
allowed bigger than one. This is both attractive and exciting as we
can implement the method much more effectively than ever before.
Numerical experiments have confirmed our theory.

Perspectively, since the inexact RQI has intimate relations with the
simplified Jacobi-Davidson method and the former is mathematically
equivalent to the latter when a Galerkin--Krylov type solver, e.g.,
the Lanczos method, is used for solving the linear systems, we can
use the convergence theory developed here for the inexact RQI to
help understand the inexact simplified JD method.
In this respect, Simoncini and Eld$\acute{e}$n~\cite{SimonciniRQI}
have proved the mathematical equivalence of the two methods,
and Freitag and Spence~\cite{freitag08} have given a further
analysis on the methods with preconditioned Lanczos solves. Based on
the equivalence, it is significant to extend our theory in this paper
to the simplified Jacobi--Davidson method with unpreconditioned and
preconditioned Lanczos inner solves. A similar relaxation on $\xi_{k+1}$
is expected.  Meanwhile, the
inexact inverse iteration is a simpler variation of the inexact RQI,
where varying $\theta_k$'s are fixed to be a constant $\sigma$,
leading to different convergence behavior. Thus, a specific analysis
is needed. It is likely to exploit the analysis approach used in
this paper to study the inexact inverse iteration. Finally, although
we have restricted to the Hermitian case, the analysis approach in
this paper may be applied to the inexact RQI for the non-Hermitian
eigenvalue problem, where the Arnoldi method is used for non-Hermitian
inner linear systems.


\end{document}